\documentclass[a4paper,10pt]{article}

\usepackage[utf8]{inputenc}
\usepackage[dvipsnames]{xcolor}
\usepackage{geometry}
\usepackage{xfrac}
\usepackage{amsfonts}
\usepackage{graphicx}
\usepackage{amsthm}
\usepackage{amsmath}
\usepackage{enumitem}
\usepackage{tikz-cd}
\usepackage{hyperref}
\usepackage{cleveref}

\graphicspath{ {./Figures/} }

\newcommand{\Grkn}{Gr_{k,n}}

\newcommand{\hiddenegonote}[1]{}

\newcommand{\tilingSet}{\mathbf{Til}}
\newcommand{\tilingSetN}{\mathbf{Til_n}}
\newcommand{\strandDiagSet}{\mathbf{Diag}}
\newcommand{\strandDiagSetN}{\mathbf{Diag_n}}

\newcommand{\id}{\mathop{}\mathopen{}\mathrm{id}}

\newtheorem{theorem}{Theorem}[section]
\newtheorem{proposition}[theorem]{Proposition}

\theoremstyle{definition}

\theoremstyle{definition}
\newtheorem{corollary}[theorem]{Corollary}

\theoremstyle{definition}
\newtheorem{remark}[theorem]{Remark}

\theoremstyle{definition}
\newtheorem{definition}[theorem]{Definition}

\theoremstyle{definition}
\newtheorem{example}[theorem]{Example}

\theoremstyle{definition}
\newtheorem{lemma}[theorem]{Lemma}

\newenvironment{acknowledgements} {\begin{abstract}}{\end{abstract}}

\title{Bicolored tilings and the Scott map}
\author{Joel Costa da Rocha}
\date{}

\begin{document}
\maketitle

\begin{abstract}
    We define bicolored tilings as a disk with a collection of smooth curves with a coloring map on the tiles that these curves delimit. Using two transformations, we define an equivalence on tilings. We then define the Scott map which creates a bijection between a generalised version of Postnikov diagrams and bicolored tilings, preserving equivalences, and mapping the geometric exchange of a diagram to an edge flip in a tiling.
\end{abstract}

\section{Introduction}
A Postnikov diagram is a combinatorial object consisting of strands drawn inside a disk satisfying a set of conditions \cite{Postnikov}. Postnikov diagrams play an important role relating to the Grassmannian, such as describing certain clusters for cluster algebras \cite{Marsh}\cite{BaurAlastairMarsh}\cite{FominZelevinsky}, or partitioning the total non-negative part of the Grassmannian of a given type \cite{Postnikov}. One particular property of a subclass of Postnikov diagrams is their bijection with triangulations of polygons, which preserves a geometric exchange within a diagram in the form of flipping a diagonal within the respective triangulation \cite{Scott}. The Scott map, which maps triangulations to Postnikov diagrams, generalises to a map on tilings of polygons \cite[2.1]{BaurMartin}. Naturally, the question arises whether such an equivalence happens with a broader class of diagrams and tilings.

In this paper, we introduce \textit{bicolored tilings} as a more general version of tilings, and further extend the Scott map (which maps classical tilings to diagrams) to our new class of tilings. We also introduce alternating path diagrams to extend Postnikov diagrams, allowing for cycles and unoriented lenses to appear. In doing so, we are able form a bijection between tilings and alternating path diagrams, as well as performing a flip (or \textit{mutation}) of any edge in a tiling such that it corresponds to a geometric exchange in the corresponding diagram.

The main result of this paper is that any Postnikov diagram (and more generally and alternating path diagram) is the image of a tiling by the Scott map.

With Postnikov diagrams being associated to Grassmann cells and cluster algebras, the result of this paper allows us to associate to each Grassmann cell a tiling, as well as describe some clusters of cluster algebras that arise from diagrams through tilings, with the mutation of a cluster variable corresponding to a flip of an edge in the corresponding tiling. Bicolored tilings also give us a different look at cluster algebras arising from disconnected Postnikov diagrams and connected diagrams behave in some ways as if they were disconnected. However, these topics exceed the scope of this paper, and in this paper we will focus only on the construction of bicolored tiling and their association with Postnikov diagrams, and we treat semi-connected diagrams in future work.

\begin{acknowledgements}
    The author would like to thank Karin Baur and João Faria Martins for their help and insightful comments during the writing of this paper and the work that preceded it. The author was supported by Royal Society Wolfson Fellowship 180004.
    
    The author would also like to thank the Isaac Newton Institute for Mathematical Sciences, Cambridge, for support and hospitality during the programme \textit{Cluster algebras and representation theory} where work on this paper was undertaken. This work was supported by EPSRC grant no EP/K032208/1.
\end{acknowledgements}

\section{Alternating path diagrams}

We first recall the definition of Postnikov diagrams. We do so by first introducing the more general version of alternating path diagrams, and defining Postnikov diagrams as path diagrams with additional restrictions.

\begin{definition}
    Consider a disk with $n$ vertices drawn on its boundary, labeled by the elements in $\{1,\dots,n\}$, in clockwise order. An \textit{(alternating) path diagram} consists of a finite collection of oriented paths, such that each path is either a closed cycle or strand with boundary vertices as endpoints, with each boundary vertex having exactly one incoming and outgoing strand, satisfying the following conditions:
    \begin{enumerate}[label=(\roman*)]
        \item
            A path does not cross itself in the interior of the disk.
        \item
            No three paths cross in one single point.
        \item
            All crossings are transversal (left figure as opposed to right figure).
            \begin{center}
                \includegraphics[scale=0.50]{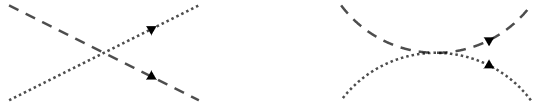}
            \end{center}
        \item
            There are finitely many crossings between paths.
        \item
            Following any path in one direction, the paths that intersect it must alternate in orientation.
            \begin{center}
                \includegraphics[scale=0.5]{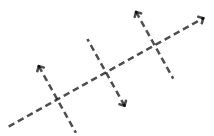}
            \end{center}
    \end{enumerate}
    We define alternating path diagrams up to equivalence of two local transformations, namely twisting and untwisting oriented lenses inside the disk or on the boundary
    \begin{center}
        \includegraphics[scale=1.5]{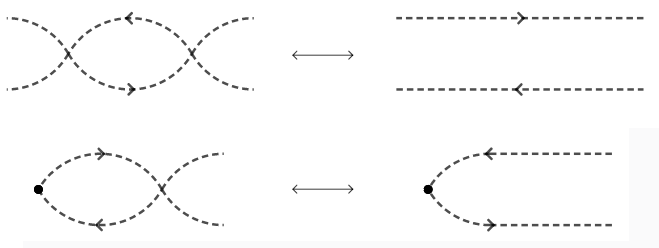}
    \end{center}
    We call the transformation from left to right a \textit{reduction}, and diagrams to which no further reduction can be applied \textit{reduced diagrams}. If two path diagrams $D_1,D_2$ are equivalent up to this transformation, we write $D_1 \equiv D_2$. The set of alternating path diagrams in a disk with $n$ boundary vertices up to equivalence is denoted $\strandDiagSetN$. The set of all alternating path diagrams is denoted $\strandDiagSet = \bigcup \strandDiagSetN$.
    Furthermore, we treat every diagram up to isotopy with the boundary vertices fixed. When necessary, we denote this equivalence $\sim$.
    For any $i \in \{1,\dots,n\}$, the strand that starts at the boundary vertex $i$ is denoted $\gamma_i$.
\end{definition}

\begin{example}
    An alternating path diagram has always $n$ strands with endpoints at $\{1,\dots,n\}$ and might also have internal cyclical paths.
    \begin{center}
        \includegraphics[scale = 1.5]{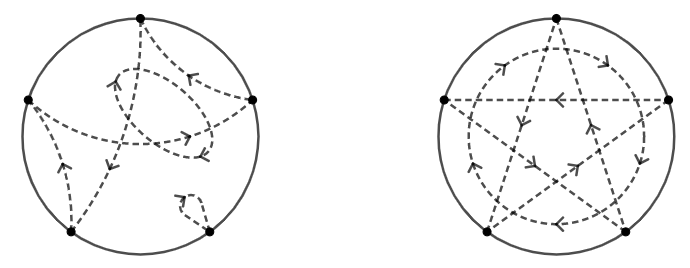}
\end{center}
\end{example}

\begin{definition}
    A \textit{Postnikov diagram}, or (alternating) strand diagram is an alternating path diagram such that
    \begin{enumerate}[label=(\roman*)]
        \item
            There are no closed cycles, i.e. every path is a strand that starts at a boundary vertex and ends in a boundary vertex. Equivalently, there are as many paths as there are boundary vertices.
        \item
            If two strands cross at two points $A$ and $B$, then one strand is oriented from $A$ to $B$, and the other from $B$ to $A$ (left figure as opposed to right figure). In other words, no two strands create unoriented lenses.
            \begin{center}
                \includegraphics[scale=0.50]{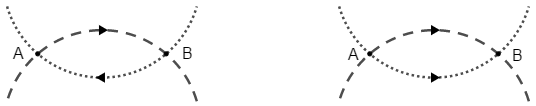}
            \end{center}
    \end{enumerate}
\end{definition}

\begin{remark}
    For the rest of this paper, when using the term \textit{diagram}, we refer to alternating path diagrams.
\end{remark}


\begin{definition} \cite{Postnikov}
    A \textit{decorated permutation} $\overline{\pi}$ of $\{1,\dots,n\}$ is a pair $(\pi,c)$ consisting of a permutation $\pi$ of $\{1,\dots,n\}$ and a map $c$ that maps any fixed point of $\pi$ to an element in $\{-1,1\}$.
\end{definition}

The function $c$ is a colouring of the fixed points of $\pi$. Any diagram defines a permutation $\pi$ of $\{1,\dots,n\}$ where $\pi(i) = j$ when $\gamma_i$ ends at the boundary vertex $j$. For any fixed point $i$ where $\gamma_i$ is oriented clockwise, $col(i) = 1$. Otherwise $col(i)=-1$. We call $i$ the \textit{source (vertex)} and $\pi(i)=j$ the \textit{target (vertex)} of $\gamma_i$.

Furthermore, any Postnikov diagram subdivides the disk into alternating and oriented regions. This is a direct consequence of the alternating property: it is easy to observe by noticing that every other strand around a common face must have the same orientation with respect to that face. For any $i \in \{1,\dots,n\}$, the alternating regions left of the strand $\gamma_i$ are given the label $i$. This will equip each region with the same number of labels \cite[9.4.1]{Marsh}. We say that the Postnikov diagram is of \textit{type} $(k,n)$ if the disk has $n$ boundary vertices and each alternating region has $k$ labels. A Postnikov diagram with permutation of type $i \longmapsto i+k \mod n$ for some fixed $k$ will always have $k$ labels per alternating region \cite[9.4.2(b)]{Marsh}, and will also be called $\Gamma_{k,n}$-diagram.

\begin{example}
The following Postnikov diagram which appeared in \cite{BaurAlastairMarsh} is a $\Gamma_{3,7}$-diagram with permutation $i \mapsto i + 3$.
    \begin{center}
        \includegraphics[scale = 0.2]{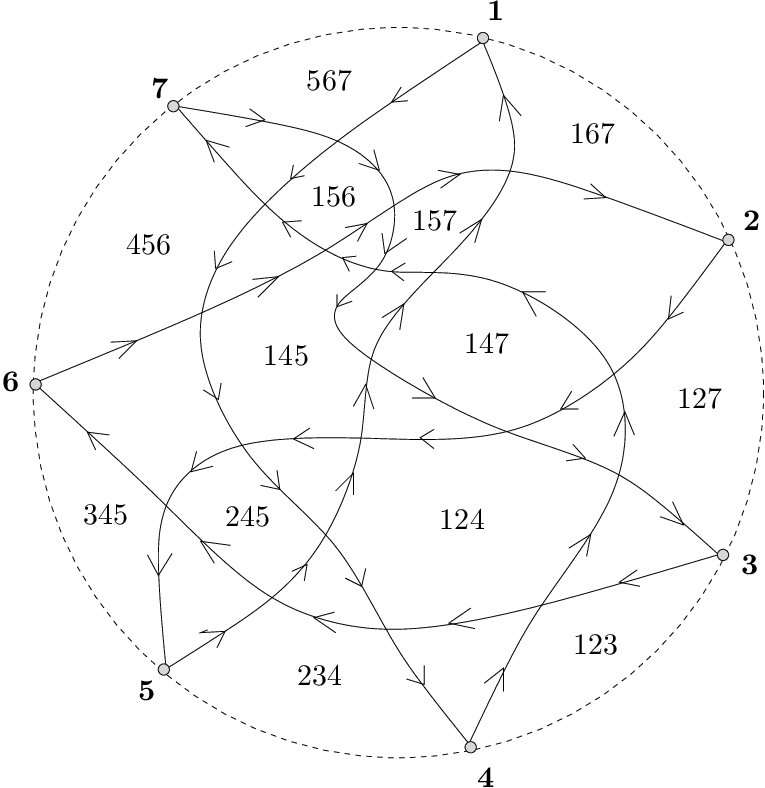}
    \end{center}
\end{example}

\begin{definition}
    Let $\Gamma$ be a path diagram, and $\Delta$ an internal quadrilateral alternating region of the diagram. We define the \textit{geometric exchange} of $\Gamma$ with respect to $\Delta$ as the local transformation on the diagram illustrated as follows
    \begin{center}
        \includegraphics[scale=0.55]{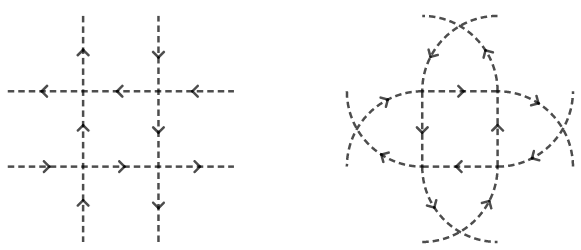}
    \end{center}
    We denote the resulting diagram $\mu_{\Delta}(\Gamma)$. It is easy to check that the geometric exchange does not change the type or permutation of a diagram. \cite{Postnikov}
\end{definition}

\hiddenegonote{
\begin{remark}
    Labels are Plücker coordinates, geometric exchange corresponds to mutation in quiver. Plücker coordinates of mutated vertices or exchanged alternating regions change with respect to Plücker relations.
\end{remark}}

\section{Bicolored tilings}

In this section, we introduce the notion of bicolored tilings. We will use them later as a source for alternating path diagrams.

\begin{definition}
    Let $D_n$ be a disk with $n$ distinct boundary vertices, enumerated $\{1,\dots,n\}$ in a clockwise order, and $x_1,\dots,x_m$ internal vertices. Then a \textit{tiling} $T$ of $D$ is $D$ equipped with a finite collection of pairwise non-crossing smooth open curves, called \textit{edges}, contained in $D$ with the endpoints of each curve (i.e. the boundaries the curves) being vertices of $D$, such that each vertex is incident to at least one edge.
    
    A tiling separates the disk into regions called \textit{tiles}, which are maximal connected subsets of the disk not intersecting any curves. Informally, they are the regions that the edges delimit in the circle. We say that a tile is of \textit{size} $s$ or an $s$-gon, if it is delimited by $s$ curves.
    
    We treat every tiling up to isotopy with the boundary vertices fixed. When necessary, we denote this isotopy $\sim$. We will consider tilings of $D = D_n$ with varying numbers of internal vertices. As we will later see, we can define an equivalence on tilings where tilings with a different number of internal vertices may be equivalent.
\end{definition}

\begin{remark}
    The $n$ boundary vertices segment the boundary of the disk into $n$ curves that we treat as part of the tiling and call \textit{boundary edges}.
\end{remark}

\begin{definition}
    A \textit{bicolored tiling} is a tiling $T$ equipped with a coloring map $c:\{\text{tiles of } T\} \longrightarrow \{0,1\}$. We think of a tile that maps to $0$ as being white and a tile that maps to $1$ as being black (or shaded dark grey).
\end{definition}

\begin{remark}
    We treat black tiles as if they were generalised versions of edges. In other words, a common edge may also be treated as a black $2$-gon and also drawn as such using curved lines. Furthermore, black tiles merge along a common edge, e.g. two black tiles of size $a$ and $b$ that are adjacent along an edge are treated as one black tile of size $a+b-2$. Similarly, a black tile of size $d$ can be split into two tiles of size $a,b$ with $a+b = d+2$.
    \begin{center}
        \includegraphics[scale=0.6]{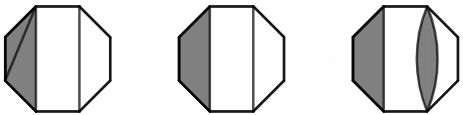}
    \end{center}
    Under merging and splitting, these tilings are all considered equal.
\end{remark}

We now introduce an equivalence class of tilings of discs with a fixed number of boundary vertices, but a varying number of internal vertices.

\begin{definition}
    Let $T$ be a bicolored tiling of $D_n$ and $p$ a tile of $T$, with edges $e_1,\dots,e_m$, for $m \geq 2$. We add two black triangular tiles $t,t'$ inside $p$ to $T$ such that
    \begin{enumerate}[label=(\roman*)]
        \item $t$,$t'$ are adjacent to distinct edges $e_i$,$e_j$ respectively.
        \item $t,t'$ share a common vertex inside the polygon, namely the vertex that is not adjacent to $e_i$ and $e_j$. 
    \end{enumerate}
    Then we call the addition or removal or such triangles an \textit{hourglass movement/transformation} on the tiling.
    \begin{center}
        \includegraphics[scale=0.65]{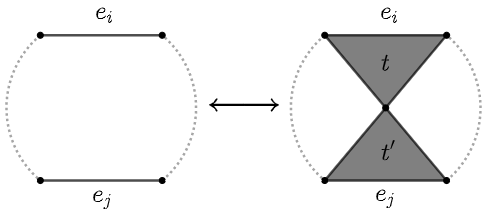}
    \end{center}
\end{definition}

\begin{definition}
    Let $T$ be a bicolored tilings and $d$ a white $2$-gon with vertices $x_1$ and $x_2$, and with black tiles $b_1$ and $b_2$ on either side. We contract the $2$-gon $d$, which shrinks the size of $b_1$ and $b_2$ by $1$ to black tiles $b_1'$ and $b_2'$, and contracts the vertices $x_1$ and $x_2$ to one single vertex $x$. Then we call this a \textit{digon contraction}, or \textit{digon decontraction} if we go the other way, of the tiling.
    \begin{center}
        \includegraphics[scale=0.65]{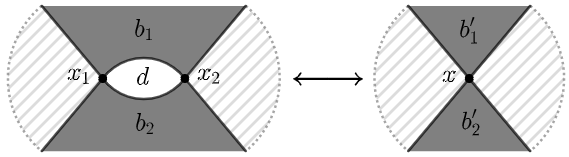}
    \end{center}
\end{definition}

Note that both the hourglass movement and the digon (de)contraction change the number of internal vertices.

\begin{definition}
    We say that two bicolored tilings $T,T'$ are \textit{equivalent}, and denote $T \equiv T'$, if $T$ is obtained from $T'$ by finitely many hourglass movements, digon contractions, and digon decontractions.
\end{definition}

The following result is straightforward.
\begin{lemma}
    The relation $\equiv$ is an equivalence relation on bicolored tilings.
\end{lemma}

\begin{definition}
    Unless otherwise specified, when talking about tilings, we mean bicolored tilings up to equivalence. The set of all bicolored tilings up to equivalence is denoted $\tilingSet$, the set of tilings with $n$ boundary vertices is denoted $\tilingSetN$.
\end{definition}

\begin{example} The following four tilings are equivalent. Note that to go from the third to the fourth tiling, the added the two black triangles merge into one bigger tile.
    \begin{center}
        \includegraphics[scale=0.65]{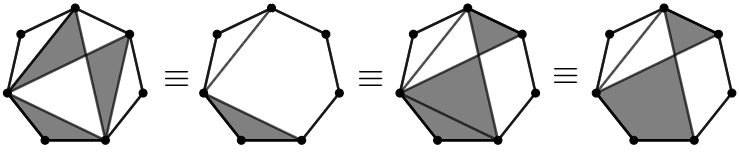}
    \end{center}
\end{example}

\begin{example} Two special cases of the hourglass equivalence are the following
    \begin{center}
        \includegraphics[scale=0.65]{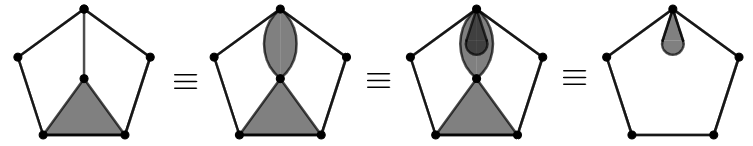}
    \end{center}
    The vertical edge can be viewed as a black $2$-gon, which can then be viewed as a black triangle (shaded in a lighter grey) added along the edge of a black $1$-gon (shaded in darker grey). With the bottom triangle, this is another hourglass.
    
    Similarly, an internal vertex with two incident edges is another hourglass that can be removed or added.
    \begin{center}
        \includegraphics[scale=0.65]{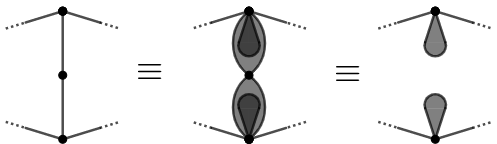}
    \end{center}
\end{example}

\begin{example}
    The following two tilings are equivalent. We can go from left to right and vice versa by contracting and decontracting the digon respectively.
    \begin{center}
        \includegraphics[scale=0.55]{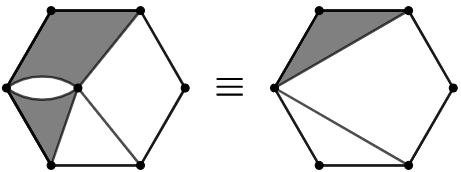}
    \end{center}
\end{example}

\section{Scott map on bicolored tilings}

We use bicolored tilings to obtain path diagrams. The Scott map defined on triangulations in \cite[p.14]{Scott} was extended to a bigger class of tilings in \cite[2.1]{BaurMartin}. We further generalise the Scott map to a map on bicolored tilings. To define the generalised map, we use a language close to the one used in \cite{Scott}.

\begin{definition}
    We define the \textit{Scott map}
    $$S: \tilingSetN \longrightarrow \strandDiagSetN, T \longmapsto D$$
    such that
    \begin{enumerate}[label=$\cdot$]
        \item any white tile is mapped to a path configuration consisting of $m$ paths, where $m$ is the size of the tile, following around the border in a counter-clockwise orientation. For example:
        \begin{center}
            \includegraphics[scale=0.54]{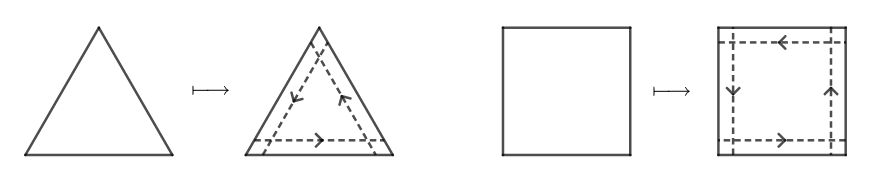}
        \end{center}
        \item any black tile is mapped to a path configuration consisting of $m$ paths, where $m$ is the number of vertices of the tile, such that each path forms an arc around a vertex inside the tile in a clockwise orientation. For example:
        \begin{center}
            \includegraphics[scale=0.54]{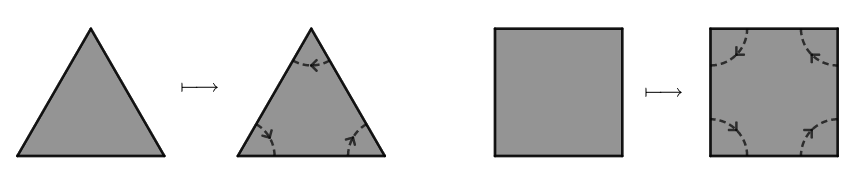}
        \end{center}
        \item If two tiles are adjacent, join the pairs of oriented paths along their shared boundary. The oriented curves obtained from concatenating all paths make up the full paths of the diagram. One can check that these are indeed consistently oriented.
        \item The paths join at the boundary, i.e. for any boundary vertex, we take the two paths that intersect the boundary on either side of the vertex closest to it and join them.
    \end{enumerate}
\end{definition}

The Scott map is invariant under merging and splitting.

\begin{example} We consider the image of a bicolored tiling of an octagon with one internal vertex under the Scott map. The image is a Postnikov diagram of type $(5,8)$.
    \begin{center}
        \includegraphics[scale=0.65]{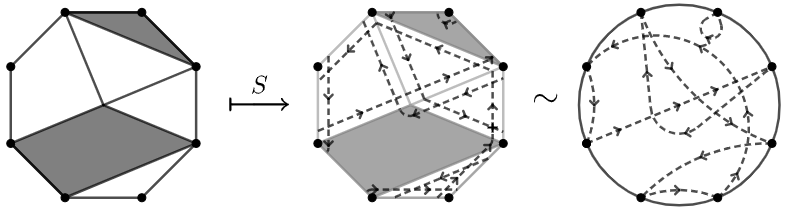}
    \end{center}
\end{example}

\begin{example}
    Not every tiling maps to a Postnikov diagram. The following tiling maps to a path diagram with two closed cycles and several unoriented lenses.
    \begin{center}
        \includegraphics[scale=1.85]{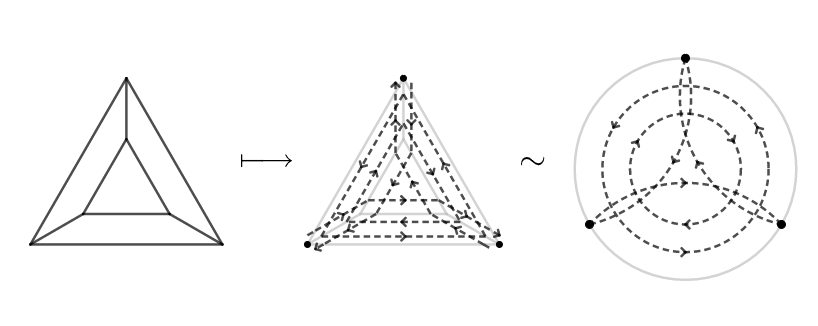}
    \end{center}
\end{example}

\begin{lemma}
    For any tiling $T$, $S(T)$ is an alternating path diagram.
    \begin{proof}
        By construction, all intersections of paths are given by the intersections of paths within white tiles. They are trivially pairwise, transversal, and there are finitely many. It is also easy to see that that the endpoints of strands are on the boundary of $T$. Any endpoints of a strand can only arise as the endpoint of a path in the configuration if that endpoint is not attached to another path. However, that is only the case if there is no tile adjacent to the boundary edge of the tile on which that endpoint lies, which means that the endpoint lies on the boundary of $T$.
        
        Following any strand from its source vertex to its target vertex, any crossings occur only in white tiles and come in pairs of two, always crossing from left to right first (after entering the tile), and from right to left second (before leaving the tile). Thus the alternating property is satisfied for any strand in $T$. If the path is a cycle, we choose any point of the cycle that lies on an edge of a white tile and treat it as both source and target vertex, and the same reasoning follows. Thus the alternating property is satisfied everywhere in $T$.
        
        Thus the resulting construction is an alternating path diagram.
    \end{proof}
\end{lemma}

\begin{remark}
    There is a correspondence between the parts that make up a tiling $T$ and the regions of the diagram $\Gamma = S(T)$. More precisely, if $v$ is a vertex, $e$ an edge or black tile, and $f$ a face in $T$, then
    \begin{enumerate} [label = $\cdot$]
        \item $v$ maps to a clockwise oriented region in $\Gamma$.
        \item $e$ maps to an alternating region in $\Gamma$.
        \item $f$ maps to a counterclockwise oriented region in $\Gamma$.
    \end{enumerate}
\end{remark}

\begin{proposition}  \label{scottPreservesEquiv}
    Let $T,T' \in \tilingSet$. If $T \equiv T'$ then $S(T) \equiv S(T')$.
    \begin{proof}
        Since both relations are equivalence relations, it is sufficient to show that if $T'$ is obtained from $T$ by the addition of an hourglass or by digon (de)contraction, then $S(T) \equiv S(T')$. Indeed, adding an hourglass results in locally twisting two paths of the diagram into a clockwise oriented lens,.
        \begin{center}
            \includegraphics[scale = 0.65]{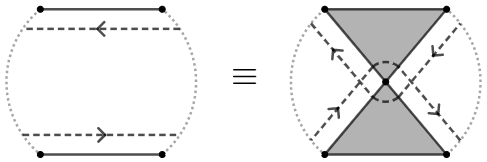}
        \end{center}
        whereas contracting a digon results in locally untwisting two paths of the diagram into a counterclockwise oriented lens.
        \begin{center}
            \includegraphics[scale=0.65]{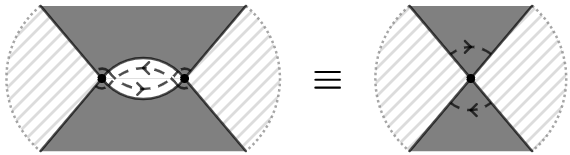}
        \end{center}
        In both cases, $S(T) \equiv S(T')$.
    \end{proof}
\end{proposition}

\begin{definition}
    We say that a tiling $T$ has \textit{decorated permutation} $\pi$ if $S(T)$ has decorated permutation $\pi$. Furthermore, we say that $T$ has \textit{type} $(k,n)$ if $S(T)$ is a Postnikov diagram of type $(k,n)$.
\end{definition}

\begin{remark}
    If $V$ is the number of all vertices and $F$ the number of white tiles of a tiling $T$ with type $(k,n)$, then $k = V-F$. This equality allows us to define the type of the tiling independently of diagrams, and extends the notion of \textit{type} to more tilings, regardless of what diagram they map to under the Scott map. The proof of this equality, however, requires introducing the notion of plabic graphs and how they parametrize positroid cells in the Grassmannian $\Grkn$ \cite{Williams}, which exceeds the scope of this paper, and thus will be treated in future work. Nevertheless, it may be helpful to keep that identity in mind.
    
    It is also worth noting that both the hourglass movement and digon (de)contraction change the number of internal vertices, but also change the number of white tiles, and that $V-F$ is invariant under both.
\end{remark}

\begin{example}
    If we arrange $n$ white $2$-gons around a common center in an $n$-gon, with the other vertex of each two $2$-gon being a boundary vertex, we obtain a tile of type $(1,n)$ and permutation $i \longmapsto i+1 \mod n$. For $n=4$, this looks as follows.
    \begin{center}
        \includegraphics[scale=0.56]{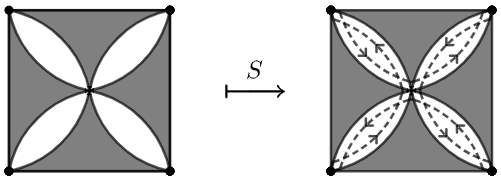}
    \end{center}
    We call this construction an \textit{$n$-antigon} (mirroring white $n$-gons whose permutation is $i \longmapsto i-1 \mod n$), and it will be relevant later when we define an inverse to the Scott map.
\end{example}

\begin{example}
    A triangulation of an $n$-gon without internal vertices is of type $(2,n)$ \cite{Scott}[Cor.2] and has permutation $i \mapsto i+2$.
    \begin{center}
        \includegraphics[scale = 0.5]{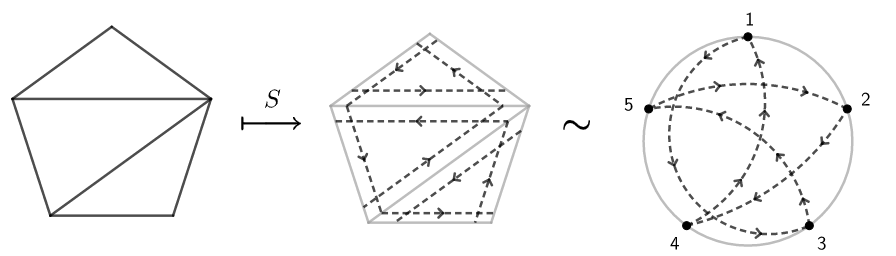}
    \end{center}
\end{example}

\begin{corollary}
    If $T \equiv T'$, then $T$ and $T'$ have the same type and permutation.
    \begin{proof}
        This directly follows by $S(T) \equiv S(T')$ (\cref{scottPreservesEquiv}), and the fact that equivalent diagrams have the same type and permutation.
    \end{proof}
\end{corollary}

\begin{definition}
    Let $T$ be a tiling. Let $e$ be an edge in $T$ that is adjacent to two triangular white tiles $t_1,t_2$. Let $v_1,v_2$ be the vertices adjacent to $e$, and $u_1,u_2$ the remaining two vertices of $t_1$ and $t_2$ respectively that are not adjacent to $e_1,e_2$. Then we define the \textit{mutation} of $T$ with respect to $e$ as the tiling $T'$ obtained by removing $e$ and adding an edge $e'$ between $u_1$ and
    $u_2$.
    \begin{center}
        \includegraphics[scale=0.60]{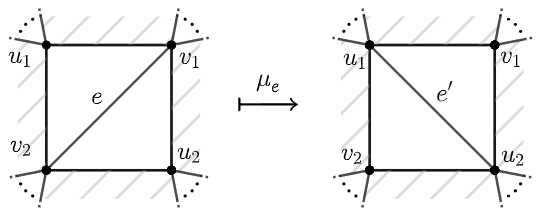}
    \end{center}
    We say that we \textit{flip} or \textit{mutate} the edge $e$ in $T$, and denote $\mu_e (T) = T'$. By convention, we may denote the new edge $e'$ as $e$. Flipping the same edge twice returns the original tilings, i.e. for any edge $e$ in $T$, $\mu_e^2(T)=\mu_e(\mu_e(T)) = T$.
    
    The setup with bicolored tilings corresponds to the classical notion of an edge flip in triangulations, we merely extend the flip to any edge inscribed in a quadrilateral in bicolored tilings. This allows us to flip edges between two white tiles of arbitrary size (in other words edges that are diagonals of white polygons within the tiling) using the hourglass equivalence, as we will see below.
\end{definition}

\begin{definition}
    Two tilings $T$ and $T'$ are said to be \textit{flip-equivalent} if one can be obtained from the other by a finite sequence of edge flips. The \textit{flip/mutation equivalence class} of $T$ is the set of all tilings that are flip-equivalent to T.
\end{definition}

\begin{example}
    The flip-equivalence class of the triangulation of a pentagon, which is of type $(2,5)$, consists of the following five tilings.
    \begin{center}
        \includegraphics[scale=0.55]{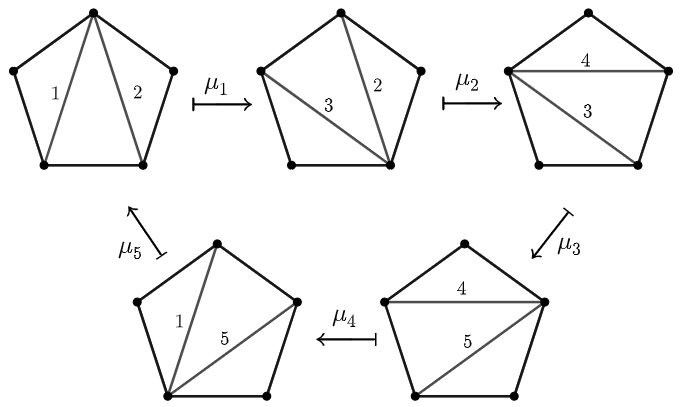}
    \end{center}
\end{example}

\begin{remark} \label{combinedRemark}
    The flip of an edge in a quadrilateral corresponds to a geometric exchange in the strand diagram that the tiling maps to, i.e. $S(\mu_e(T)) \equiv \mu_{S(e)}S(T)$ for any edge $e$ in a tiling $T$. This is pointed out in \cite{Scott} and still remains true in bicolored tilings.
    \begin{center}
        \includegraphics[scale=0.45]{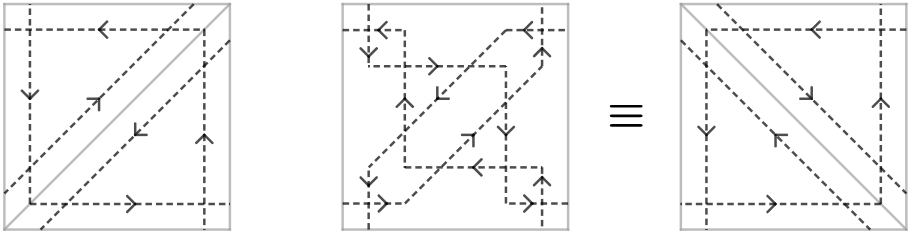}
    \end{center}
    Postnikov diagrams naturally map to quivers, which generate cluster algebras, with the vertices of the quiver corresponding to cluster variables \cite{BaurAlastairMarsh}. The geometric exchange and consequently the flip of an edge correspond to the mutations of cluster variables, with each tiling corresponding to a cluster.
\end{remark}

As noted before, we can use the hourglass equivalence of tilings to flip edges between any two white tiles of size greater than $2$. We may do that by adding hourglasses around the edge that we want to flip in order to inscribe it in a white quadrilateral, which may then be flipped.
\begin{center}
    \includegraphics[scale = 0.55]{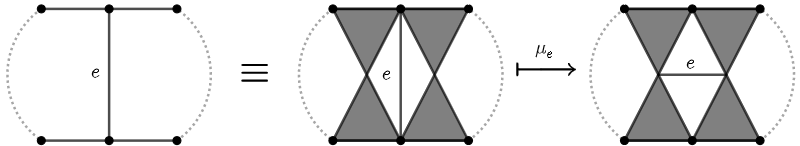}
\end{center}
In other words, any edge that acts as the diagonal of a white polygon can be flipped that way.

\begin{example}
    We consider the following rhombic tiling with internal edges labeled $e_1,e_2,e_3$ clockwise, starting with the left one. To mutate the edge $e_1$, we add two hourglasses inside the adjacent tiles, which then circumscribe the edge inside a quadrilateral, allowing us to flip it. We denote $\mu_i := \mu_{e_i}$.
    
    \begin{center}
        \includegraphics[scale = 0.65]{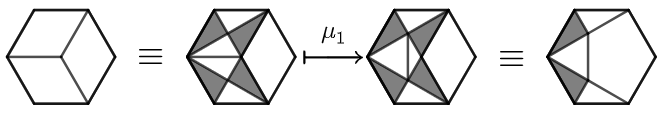}
    \end{center}
    
    Further flipping the edges $e_2,e_3,e_1$ in that order indeed returns another rhombic tiling, which is the original tiling to which we applied a Yang-Baxter move.
    \begin{center}
        \includegraphics[scale = 0.65]{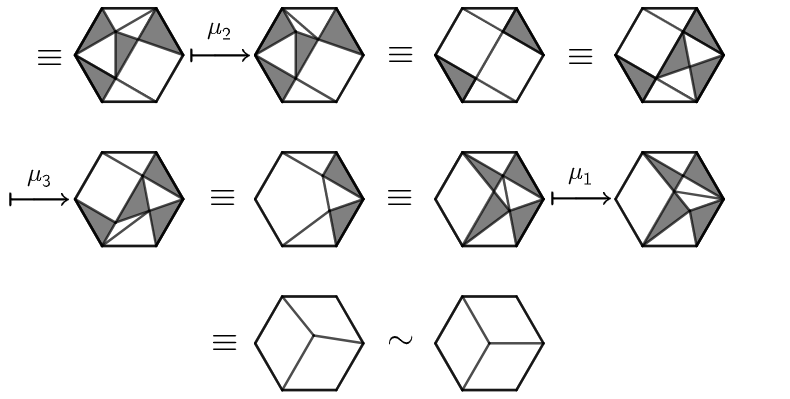}
    \end{center}
\end{example}

\section{Inverse Scott map}

By defining an inverse to the Scott map, we can, in particular, construct a tiling for any Postnikov diagram, and as a consequence construct a tiling of any type and permutation as well.

\begin{definition}
    Let $D \in \strandDiagSetN$ be an alternating path diagram, and let $r$ be a region of $D$, then the \textit{size} of $r$ is the number of crossings in $D$ which lie on the boundary of $r$.
\end{definition}

\begin{definition}
    We define the map
    $$\overline{S}: \strandDiagSetN \longrightarrow \tilingSetN$$
    that maps any path diagram to a tiling of an $n$-gon such that
    \begin{enumerate}[label = $\cdot$]
        \item any counterclockwise region of size $m$ is mapped to a white $m$-gon.
        \item any clockwise region of size $m$ is mapped to an $m$-antigon for $m > 1$, while a clockwise region of size $1$ (i.e. a clockwise loop) is mapped to a black $1$-gon.
        \item any alternating and any boundary region of size $m$ is mapped to a black $m$-gon.
        \item if two regions are edge- or vertex-adjacent, then their images are also edge- or vertex-adjacent, respectively.
    \end{enumerate}
\end{definition}

\begin{proposition}
    $\overline{S}$ is the inverse of $S$, i.e. $S \overline{S} = \id_{\strandDiagSetN}$ and $\overline{S}S = \id_{\tilingSetN}$. We call $S^{-1} = \overline{S}$ the \textit{inverse Scott map}.
    \begin{proof}
        Since both maps are defined locally, we can simply look at how applying those maps in succession locally returns the initial construction.
        \begin{enumerate}[label = $\cdot$]
            \item $S\overline{S} = \id_{\strandDiagSetN}$: We look at the intersection of two paths. The paths divide the region locally into 4 smaller regions, one of which is clockwise oriented. Assume that the clockwise-oriented region is of size $m>1$, which means that $\overline{S}$ maps it to an $m$-antigon. We see that applying $\overline{S}$ followed by $S$ returns the same intersecting sections of paths 
            \begin{center}
                \includegraphics[scale=0.60]{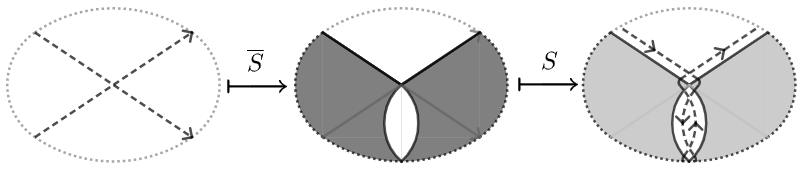}
                \includegraphics[scale=0.60]{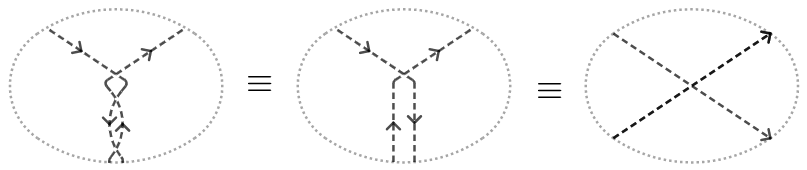}
            \end{center}
            Thus $S\overline{S}$ does not locally change the diagram, and therefore maps the diagram to itself. The case where $m=1$ is similar.
            \item $\overline{S}S = \id_{\tilingSetN}$: Similarly, we look at how $\overline{S}S$ applies locally on a tiling. More precisely, we look at where intersections of the diagram would occur after mapping the tiling by $S$, which would be at the corner of each white tile. In the illustration below the white tile is adjacent to two black tiles around that corner. However, as seen before, edges can be treated as black $2$-gons and drawn as such, hence we do not lose generality with this construction.
            \begin{center}
                \includegraphics[scale=0.65]{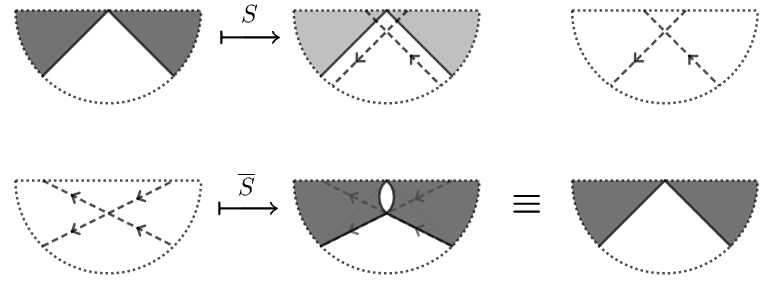}
            \end{center}
            We observe that $\overline{S}S$ does not locally change the diagram, and therefore maps the diagram to itself.
            
            Once again, in the left figure of second line, we assume that $m>1$, where $m$ is the size of the clockwise oriented region, which maps the region to an antigon under $\overline{S}$. The result is similar if $m=1$.
        \end{enumerate}
        Thus we conclude that $\overline{S} = S^{-1}$.
    \end{proof}
\end{proposition}

\hiddenegonote{Prove that the is well-defined. The definition itself does not guarantee that there is one unique tiling (up to equivalence) that satisfies the 4th condition, it could be possible that you could arrange the images of the regions in a way to satisfy the adjacencies in two different ways.}

\begin{example} \label{diag36} We show the effect of $S^{-1}$ on a diagram of type $(3,6)$
    \begin{center}
        \includegraphics[scale=1.2]{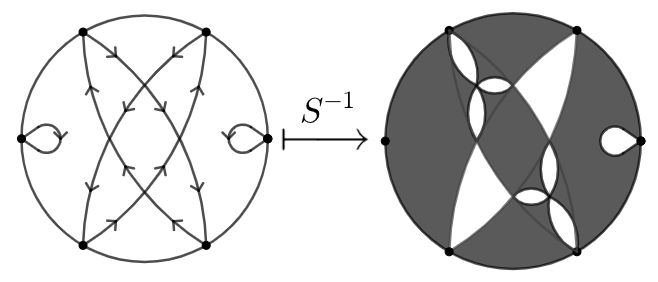}
    \end{center}
\end{example}

\begin{remark}
    The inverse Scott map as defined above helps us determine that every diagram has a corresponding tiling, i.e. a tiling that maps to that diagram. However, the result is not always reduced, and we may have to use the equivalence on tilings to obtain a "cleaner" result.
\end{remark}

\begin{example} Using the equivalence on tilings, we obtain the following tiling for the strand diagram of \cref{diag36}.
    \begin{center}
        \includegraphics[scale=1.0]{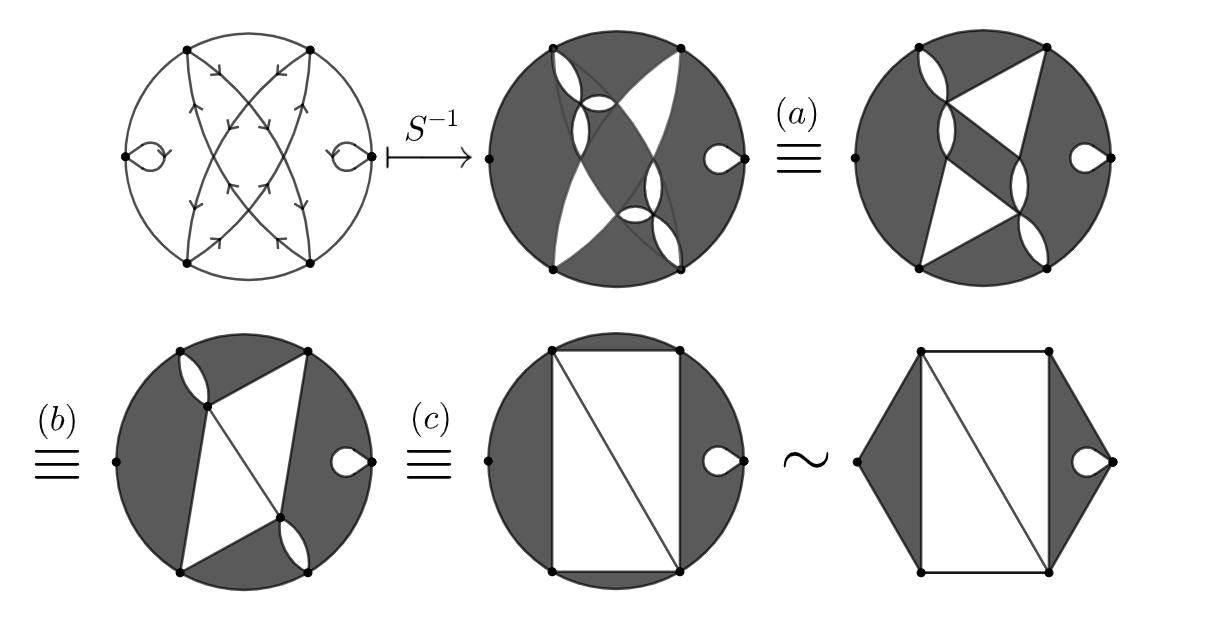}
    \end{center}

\begin{enumerate}[label = (\alph*)]
    \item We contract the two horizontal lenses to the left/right of the top-right/bottom-left white triangle, respectively.
    \item We contract the two vertical lenses to the bottom/top of the top-right/bottom-left triangle, respectively.
    \item We contract the two remaining lenses on the boundary.
\end{enumerate}
\end{example}

\begin{theorem}
    For any decorated permutation $\pi$, there is a tiling $T$ whose permutation is $\pi$. Alternatively, any Postnikov diagram can be obtained as the image of a bicolored tiling by the Scott map.
    \begin{proof}
        For any decorated permutation $\pi$, there is a Postnikov diagram $\Gamma$ with permutation $\pi$ \cite{Postnikov}[14.2, 14.4 -14.7]. If $T = S^{-1}(\Gamma)$, then $S(T) = S(S^{-1}(\Gamma)) = \Gamma$. And thus $T$ maps to $\Gamma$ by the Scott map and, by definition, $T$ has permutation $\pi$.
    \end{proof}
\end{theorem}

\section{Final remarks}

In this paper we focused on introducing the class of bicolored tilings and constructing a correspondence with alternating strand diagrams and decorated permutations. Some applications and properties of this construction will be explored in other papers, as certain obstacles remain with bicolored tilings and the generalised Scott map. For instance, the Scott map does not guarantee that a tiling maps to a Postnikov diagram, which would motivate describing which tilings map to Postnikov diagrams in particular. Another point of interest would be to find a method to generate tilings for any given decorated permutation, or whether certain permutations can be obtained by some easily identifiable classes of tilings.

\newpage

\bibliographystyle{plain}
\bibliography{main}

\end{document}